\documentclass[11pt]{article}

\usepackage{amscd,amsmath, amssymb, fancyhdr}

\numberwithin{equation}{section}

\newcommand{\version}{version 1.0,\ \ August 22, 2019}

\def\eqref#1{(\ref{#1})}

\newcommand{\arrow}{{\:\longrightarrow\:}}

\def\C{{\Bbb C}}
\newcommand{\R}{{\Bbb R}}

\def\1{\sqrt{-1}\:}

\newcommand{\cntrct}                
{\hspace{2pt}\raisebox{1pt}{\text{$\lrcorner$}}\hspace{2pt}}

\renewcommand{\bar}{\overline}
\renewcommand{\phi}{\varphi}
\renewcommand{\epsilon}{\varepsilon}
\renewcommand{\geq}{\geqslant}
\renewcommand{\leq}{\leqslant}

\newcommand{\Def}{\operatorname{Def}}

\newcounter{Mycounter}[section]
\newcounter{lemma}[section]
\renewcommand{\thelemma}{{Lemma \thesection.\arabic{lemma}}}
\newcommand{\lemma}{%
    \setcounter{lemma}{\value{Mycounter}}
    \refstepcounter{lemma}
    \stepcounter{Mycounter}
    {\noindent \bf \thelemma:\ }}

\newcounter{claim}[section]
\newcounter{sublemma}[section]
\newcounter{corollary}[section]
\newcounter{theorem}[section]
\renewcommand{\thetheorem}{{Theorem \thesection.\arabic{theorem}}}
\newcommand{\theorem}{%
    \setcounter{theorem}{\value{Mycounter}}
    \refstepcounter{theorem}
    \stepcounter{Mycounter}
    {\noindent \bf \thetheorem:\ }}

\newcounter{conjecture}[section]
\newcounter{proposition}[section]
\newcounter{definition}[section]
\renewcommand{\thedefinition}
      {{Definition~\thesection.\arabic{definition}}}
\newcommand{\definition}{%
    \setcounter{definition}{\value{Mycounter}}
    \refstepcounter{definition}
    \stepcounter{Mycounter}
    {\noindent \bf \thedefinition:\ }}

\newcounter{example}[section]
\renewcommand{\theexample}{{Example \thesection.\arabic{example}}}
\newcommand{\example}{%
    \setcounter{example}{\value{Mycounter}}
    \refstepcounter{example}
    \stepcounter{Mycounter}
    {\noindent \bf \theexample:\ }}

\newcounter{remark}[section]
\renewcommand{\theremark}{{Remark \thesection.\arabic{remark}}}
\newcommand{\remark}{%
    \setcounter{remark}{\value{Mycounter}}
    \refstepcounter{remark}
    \stepcounter{Mycounter}
    {\noindent \bf \theremark:\ }}

\newcounter{problem}[section]
\newcounter{question}[section]
\makeatletter

\@addtoreset{equation}{section} \@addtoreset{footnote}{section}
\makeatother

\def\blacksquare{\hbox{\vrule width 5pt height 5pt depth 0pt}}
\def\endproof{\blacksquare}

\begin{document}
\begin{center}
{\LARGE\bf
Kobayashi non-hyperbolicity of Calabi-Yau manifolds via mirror 
symmetry \\[4mm]
}

Ljudmila Kamenova\footnote{Partially supported 
by a grant from the Simons Foundation/SFARI (522730, LK).}, 
Cumrun Vafa\footnote{The research of CV is supported in part by the 
NSF grant PHY-1067976 and by a grant from the Simons Foundation (602883, CV).}

\end{center}

{\small \hspace{0.10\linewidth}
\begin{minipage}[t]{0.85\linewidth}
{\bf Abstract} \\
A compact complex manifold is Kobayashi non-hyperbolic if there exists an  
entire curve on it. Using mirror symmetry we establish that there are 
(possibly singular) elliptic or rational curves on any Calabi-Yau manifold 
$X$, whose mirror dual $\check X$ exists and is not ``Hodge degenerate'', 
therefore proving that $X$ is Kobayashi non-hyperbolic. We are not aware 
of any higher dimensional simply connected Calabi-Yau manifolds that satisfy  
the ``Hodge degenerate'' condition. 
\end{minipage}
}

\tableofcontents


\section{Introduction}


In \cite{_Kobayashi:1976_} S. Kobayashi conjectured that all compact 
Calabi-Yau manifolds have vanishing Kobayashi pseudometric. 
Kamenova-Lu-Verbitsky (\cite{klv}) proved Kobayashi's conjecture 
for all K3 surfaces and for certain compact hyperk\"ahler manifolds 
that are deformation equivalent to Lagrangian fibrations. According to the 
hyperk\"ahler SYZ conjecture, any compact hyperk\"ahler manifold is 
deformation equivalent to a Lagrangian fibration. 

\hfill

A weaker version of Kobayashi's conjecture states that any compact 
Calabi-Yau manifold $X$ is Kobayashi non-hyperbolic, i.e., there is a 
non-constant holomorphic map from $\C$ to $X$. This would be the case 
if there are (possibly singular) rational or elliptic curves on $X$. 
Using ergodicity of hyperk\"ahler complex structures, in 
\cite{_Verbitsky:ergodic_}
M. Verbitsky showed that all compact hyperk\"ahler manifolds are 
Kobayashi non-hyperbolic. Using mirror symmetry 
(more precisely, relating Ray-Singer torsion and genus one 
Gromov-Witten invarians) we show that a Calabi-Yau manifold $X$ is Kobayashi 
non-hyperbolic, assuming its mirror dual $\check X$ exists and satisfies 
some non-degeneracy condition involving Hodge 
bundles on the deformation space $\Def(\check{X})$ of $\check X$. 
Other than K3 surfaces and complex tori, we are not aware of other Calabi-Yau 
manifolds that satisfy the ``Hodge degenerate'' condition. Even though 
K3 surfaces and complex tori are ``Hodge degenerate'', they are 
Kobayashi non-hyperbolic. We suspect that, as Kobayashi conjectured,  
all Calabi-Yau manifolds are Kobayashi non-hyperbolic and that the 
``Hodge degenerate'' condition is sufficient but not necessary. 

\hfill

In \cite{hbw} D. R. Heath-Brown and P. Wilson have shown that for a Calabi-Yau 
threefold $X$ with a large Picard number $\rho(X) > 13$ there is a birational 
morphism $\phi: X \arrow \bar X$ onto a normal projective threefold $\bar X$ 
with $\rho(\bar X) < \rho (X)$ and such that the exceptional locus of $\phi$ 
is covered by rational curves. This implies the existence of 
rational curves on $X$. In other words, $X$ is Kobayashi non-hyperbolic 
in the case of $\rho(X) > 13$. However, in the cases with low Picard number 
there wouldn't be necessarily any smooth curves, 
but there could still be singular ones. 
Using mirror symmetry one can take into account singular elliptic or 
rational curves which is sufficient to imply non-hyperbolicity under the 
additional ``Hodge non-degeneracy'' condition.


\section{Holomorphic Ray-Singer torsion and GW invariants}


In this paper we consider compact Calabi-Yau manifolds. A {\it Calabi-Yau 
manifold} $X$ is a compact complex K\"ahler manifold with a trival canonical 
bundle ${\cal K}_X$ and a trivial fundamental group $\pi_1(X)=0$. 
A {\it hyperk\"ahler manifold} is a Calabi-Yau manifold with an 
everywhere non-degenerate holomorphic $2$-form. 

\hfill

Let $X$ be a complex manifold. Recall that a {\it pseudometric} on $X$ 
is a function $d$ on $X\times X$ with values in $\R_{\geq 0}$ 
that satisfies all the properties of a 
distance function except for the non-degeneracy condition: $d(x,y)=0$ 
only if $x=y$. The {\it Kobayashi pseudometric} $d_X$ on $X$ 
is defined as the supremum of all pseudometrics $d$ on $X$ that 
satisfy the distance decreasing property with respect to holomorphic 
maps $f$ from the Poincar\'e disk $(\mathbb D,\rho)$ to 
$X$: 
\[ f^*d\leq \rho \ \mbox{, or equivalently,} \ d(f(x),f(y))\leq \rho(x,y)\ 
\, \forall x,y\in {\mathbb D}.
\]
Here $\rho$ denotes the  Poincar\'e metric on ${\Bbb D}$.

\hfill

In \cite{_Kobayashi:1976_} 
Kobayashi conjectured that if $X$ is a Calabi-Yau manifold, then 
$d_X \equiv 0$. Kamenova-Lu-Verbitsky (\cite{klv}) proved this conjecture 
for all K3 surfaces and for certain compact hyperk\"ahler manifolds 
deformation equivalent to Lagrangian fibrations. 

\hfill

\definition
A compact manifold $X$ is called {\it Kobayashi hyperbolic} 
if any holomorphic map $\C \arrow X$ is constant, or equivalently, 
the Kobayashi pseudometric is non-degenerate. A non-constant holomorphic 
map $\C \arrow X$ is called an {\it entire curve}. 

\hfill

If the Kobayashi pseudometric is degenerate, the manifold $X$ is Kobayashi 
non-hyperbolic, i.e., there are entire curves on $X$. Examples of entire 
curves include (singular) rational or elliptic curves. 
Verbitsky has shown that all compact hyperk\"ahler manifolds 
are Kobayashi non-hyperbolic, \cite{_Verbitsky:ergodic_}. 

\hfill

\definition
The {\it holomorphic Ray-Singer torsion} of a Calabi-Yau n-fold is 
\begin{equation}
T = \prod\limits_{1 \leq p,q \leq n} (\det \Delta_{p,q}')^{(-1)^{p+q}pq},
\end{equation}
where $\Delta_{p,q}'$ is the non-singular part of the $\bar{\partial}$-Laplace 
operator on the $(p,q)$-forms. 

\hfill

In the literature the holomorphic Ray-Singer torsion is also known as the 
{\it BCOV torsion}. 
Fang-Lu-Yoshikawa computed the BCOV invariant for 
quintic threefolds in $\C P^4$ explicitly (see Theorem 1.1 in \cite{fly}). 

\hfill

Let us recall the following formula from \cite{bcov} and \cite{fl}:

\hfill

\theorem
Let $\omega_{WP}$ and $\omega_{H^i}$ be the K\"ahler forms of the 
Weil-Peterson metric and the generalized Hodge metrics, respectively. 
Then 
\begin{equation} \label{ddT_formula}
\sum\limits_{i=1}^{n} (-1)^i \omega_{H^i} - \frac{\sqrt{-1}}{2\pi} 
\partial\bar\partial \log T = \frac{\chi(X)}{12} \omega_{WP},
\end{equation}
where $\chi(X)$ is the Euler characteristic of the Calabi-Yau manifold $X$. 

\hfill

Following Section 2 of \cite{fl} we define the terms in formula 
(\ref{ddT_formula}). For a polarized Calabi-Yau manifold $X$ there exists
the universal family $\frak X$ together with a proper surjective flat 
morphism $\pi: {\frak X} \arrow \Def(X)$. Let $\Omega$ be a non-zero 
holomorphic $(n,0)$-form on $X$, where $n = \dim_\C X$. 
Consider the Hodge bundles 
$PR^q\pi_{*}\Omega_{\frak X/\Def}^p\rightarrow \Def(X)$, i.e., 
the relative version of the primitive cohomology groups $PH^{p,q}(X)$. 
Then for $1 \leq i \leq n$ we have (by Proposition 2.8 in \cite{fl}) 
\begin{equation}
\omega _{H^{i}}=\sum_{1~ \leq p ~\leq i} p ~ c_{1}(R^{i-p}\pi_{*}\Omega^{p}_{\frak X/ \Def}).
\end{equation}
Similarly, the form of the Weil-Peterson metric is 
\begin{equation}
\omega_{WP}=c_1(R^0\pi_*(\Omega^n_{\mathfrak X/\Def})).
\end{equation}

\hfill

\remark
By the formulas above, 
the condition $\partial\bar\partial \log T = 0$ can be expressed in 
terms of the first Chern classes of the Hodge bundles as 
\begin{equation} \label{hdg_deg}
\sum\limits_{i=1}^{n} \sum\limits_{p=1}^{i} (-1)^i p ~ c_{1}(R^{i-p}\pi_{*} 
\Omega^{p}_{\frak X/ \Def}) = 
\frac{\chi(X)}{12} c_1(R^0\pi_*(\Omega^n_{\mathfrak X/\Def})).
\end{equation}

\hfill

\definition
We call a Calabi-Yau manifold $X$ {\em Hodge degenerate} if it satisfies 
the condition (\ref{hdg_deg}) above, relating Chern classes of Hodge bundles. 

\hfill

\remark \label{fl-small-dim}
In small dimensions there are simpler relations between the K\"ahler forms of 
the Hodge metric $\omega_H$ and the Weil-Peterson metric $\omega_{WP}$, 
namely if: 
\begin{enumerate}
\item $n=2$, then $\omega_H=2\omega_{WP}$;
\item $n=3$, then $\omega_H=(m+3)\omega_{WP}
+{\rm Ric}(\omega_{WP})$;
\item $n=4$, then $\omega_H=(2m+4)\omega_{WP}
+2{\rm Ric}(\omega_{WP})$,
\end{enumerate}
where $n=\dim X$ and $m=\dim \Def(X)$ (see \cite{bcov} and \cite{fl}). 

\hfill 

\remark 
Let $X$ be a Calabi-Yau $3$-fold. Then the action of $T(\Def(X))$ on 
the lower degree Hodge bundles is trivial, hence $\omega_{H^i}=0$ for 
all $i<3$ (see \cite{fl}). Therefore, the condition (\ref{hdg_deg}) 
becomes $-\omega_H= \frac{\chi(X)}{12} \omega_{WP}$. When we combine it 
with the second case of \ref{fl-small-dim}, we obtain the following result. 

\hfill

\lemma \label{Hodge-deg-3-fold}
A Calabi-Yau $3$-fold $X$ is Hodge degenerate if and only if 
${\rm Ric}(\omega_{WP}) = - \big( m+3+ \frac{\chi(X)}{12} \big) \omega_{WP}$, 
where $m=\dim \Def(X) = h^1(X, T_X)$.  

\hfill

\definition 
Let $\bar{\frak M}_{g,k}(X,d)$ denote the moduli space of stable degree $d$ 
maps from genus $g$ curves with $k$ marked points into a smooth projective 
variety $X$. Denote by $[\bar{\frak M}_{g,k}(X,d)]^{vir}$ the virtual 
fundamental class of $\bar{\frak M}_{g,k}(X,d)$. 

\hfill

\example \label{az}
If $X$ is a quintic threefold, then the virtual fundamental class 
$[\bar{\frak M}_{g,k}(X,d)]^{vir}$ is $0$-dimensional and its degree is 
$$N_{g,d} = \langle 1, [\bar{\frak M}_{g,k}(X,d)]^{vir} \rangle,$$ 
which is called {\it the degree $d$ genus $g$ Gromov-Witten invariant}. 
In the genus $1$ case, in $\bar{\frak M}_{1,k}(X,d)$ one can define the 
{\it main component}  $\bar{\frak M}_{1,k}^0(X,d)$ as the closure of the 
locus in $\bar{\frak M}_{1,k}(X,d)$ consisting of maps from smooth domains. 
The degree of the main component is the {\it reduced} genus 1 degree $d$ 
GW-invariant $N_{1,d}^0$. For a quintic threefold $X$, A. Zinger has established 
the following relation: $$N_{1,d} = N_{1,d}^0 + \frac{1}{12} N_{0,d}$$ (Theorem 
$1.1$ in \cite{zinger1}). In particular, since $N_{1,d} \not= 0$, then at least 
one of $N_{1,d}^0$ or $N_{0,d}$ is non-zero, and therefore there will be genus 
one or genus zero curves on $X$, which implies that $X$ is Kobayashi 
non-hyperbolic. 

\hfill

\remark
Let $N_{1,d}$ be the degree $d$ genus one Gromov-Witten invariant of a 
quintic threefold $X$. A. Zinger has computed the generating function 
$\sum\limits_{d=1}^\infty N_{1,d} e^{dt}$ on the A-side of $X$ (Theorem $1$ in 
\cite{zinger}). His computation coincides with Fang-Lu-Yoshikawa's 
computation of the BCOV torsion of the mirror dual $\check{X}$ of $X$ on 
the B-side, thus establishing that for quintic threefolds the number of 
elliptic curves on $X$ is related to the BCOV torsion of $\check{X}$ as 
predicted in Bershadsky-Cecotti-Ooguri-Vafa's paper \cite{bcov}.  

\hfill


\section{Non-hyperbolicity results}







Here we establish Kobayashi non-hyperbolicity under the additional 
Hodge non-degeneracy assumption, using the tools of mirror symmetry. 
In this section we consider Calabi-Yau manifolds with maximal holonomy $SU(n)$. 
Recall that we defined $m:=\dim \Def(X) = h^1(X, T_X)$. 

\hfill

\theorem 
Let $X$ be a Calabi-Yau $3$-fold whose mirror dual $\check{X}$ exists 
and is not Hodge degenerate, i.e., 
${\rm Ric}(\omega_{WP}) \not= - \big( m+3- \frac{\chi(X)}{12} \big) \omega_{WP}$. 
Then $X$ is Kobayashi non-hyperbolic. 

\hfill

{\bf Proof.}
Let $F_1$ be the topological string partition function of genus 
one\footnote{For example, in the case of a quintic threefold $F_1$ is 
the generating function with coefficients $N_{1,d}$, the degree $d$ genus one 
Gromov-Witten invariant of $X$.} on $X$, i.e., the topological string 
amplitude. By \cite{bcov}, $F_1$ corresponds to $\log T$, where $T$ is the 
holomorphic Ray-Singer torsion of the mirror dual $\check{X}$ of $X$. 
The contribution from constant maps in $F_1$ is 
$\frac{1}{24} \int k \wedge c_2$, where $k$ is the K\"ahler class. 
In order for $X$ to be Kobayashi non-hyperbolic (or more generally, to have 
non-constant holomorphic embeddings of genus $0$ or $1$ curves), we need that 
$F_1 - \frac{1}{24} \int k \wedge c_2 \not= 0$, for which it is enough 
that $\partial\bar\partial \log T \not= 0$ (see \cite{bcov}). 
If $\partial\bar\partial \log T \not= 0$, i.e., if $\check{X}$ is not 
Hodge degenerate, then there are rational or elliptic curves on $X$, i.e., 
$X$ is Kobayashi non-hyperbolic. Notice that since we are computing 
$\partial\bar\partial \log T$ for the mirror dual $\check{X}$ 
and $\chi (\check{X}) = -\chi(X)$, then $\check{X}$ is not Hodge 
degenerate if and only if ${\rm Ric}(\omega_{WP}) \not= - 
\big( m+3- \frac{\chi(X)}{12} \big) \omega_{WP}$ by \ref{Hodge-deg-3-fold}.
\endproof

\hfill

Similarly, if $X$ is a Calabi-Yau $n$-fold, the contribution from 
constant maps into $F_1$ is $\frac{(-1)^{n-1}}{24} \int k \wedge c_{n-1}$,  
where $k$ is the K\"ahler class, and for the condition 
$$F_1 - \frac{(-1)^{n-1}}{24} \int k \wedge c_{n-1} \not= 0$$ 
to be satisfied, it is enough to check that 
$\partial\bar\partial \log T \not= 0$, where $T$ is the holomorphic 
Ray-Singer torsion of the mirror dual $\check{X}$ of $X$.
Therefore, in the general case we obtain the same result, however 
the condition that $\check{X}$ be Hodge non-degenerate in general is 
expressed in a more cumbersome way.  

\hfill 

\theorem 
Let $X$ be a Calabi-Yau $n$-fold ($n>2$) whose mirror dual $\check{X}$ 
exists and is not Hodge degenerate, i.e., it does not satisfy equation 
(\ref{hdg_deg}). Then $X$ is Kobayashi non-hyperbolic. 

\hfill

\remark
Notice that all K3 surfaces are Hodge degenerate, because for a K3 surface $S$ 
we have $ \frac{\sqrt{-1}}{2\pi} \partial\bar\partial \log T = 
\omega_H - \frac{\chi(S)}{12} \omega_{WP} = 2 \omega_{WP} - 
\frac{24}{12} \omega_{WP} =0$ using the first case of \ref{fl-small-dim}. 
However, K3 surfaces are Kobayashi non-hyperbolic as shown in 
\cite{_Verbitsky:ergodic_}, therefore the Hodge non-degenerate condition 
is sufficient, but not necessary. In higher dimensions ($n>2$) we are not aware 
of any examples of simply connected Calabi-Yau manifolds with maximal holonomy 
that are Hodge degenerate. Complex tori, on the other hand, are Hodge 
degenerate, however they contain lots of entire curves, and therefore are 
Kobayashi non-hyperbolic. 

\hfill

{\bf Acknowledgements:} The first named author would like to thank Aleksey 
Zinger and Zhiqin Lu for the interesting discussions and references. 
We are grateful to Albrecht Klemm and Dave Morrison for their comments and 
suggestions. We also thank Martin Ro\v cek and the SCGP for their 
hospitality during the Simons Summer Workshops.

\hfill

\noindent {\sc Ljudmila Kamenova\\
Department of Mathematics, 3-115 \\
Stony Brook University \\
Stony Brook, NY 11794-3651, USA,} \\
\tt kamenova@math.stonybrook.edu
\\
 
\noindent {\sc Cumrun Vafa\\
Jefferson Physical Laboratory \\
Harvard University \\
Cambridge, MA 02138, USA,} \\
\tt vafa@physics.harvard.edu


\begin{thebibliography}{GMP}

\bibitem[BCOV]{bcov} Bershadsky, M., Cecotti, S., Ooguri, H., Vafa, C., 
Appendix by S. Katz, 
{\em Holomorphic Anomalies in Topological Field Theories,} 
Nucl. Phys. B {\bf 405} (1993) 279-304. 

\bibitem[CDFLL]{cdfll} Ceresole, A., D'Auria, R., Ferrara, S., Lerche, W., 
Louis, J., {\em Picard-Fuchs Equations and Special Geometry,} 
Int. J. Mod. Phys. {\bf A8} (1993) 79-114. 

\bibitem[FL]{fl} Fang, H., Lu, Z., 
{\em Generalized Hodge metrics and BCOV torsion on Calabi-Yau moduli,} 
J. reine angew. Math {\bf 588} (2005) 49-69. 

\bibitem[FLY]{fly} Fang, H., Lu, Z., Yoshikawa, K.-I., 
{\em Analytic torsion for Calabi-Yau threefolds,} 
J. Diff. Geom. {\bf 80} (2008) 175-259. 

\bibitem[HBW]{hbw} Heath-Brown, D. R., Wilson, P., 
{\em Calabi-Yau threefolds with $\rho > 13$,} 
Math. Ann. {\bf 294} (1992) 49-57. 

\bibitem[KLV]{klv}
Kamenova, L., Lu, S., Verbitsky, M., 
{\em Kobayashi pseudometric on hyperkahler manifolds,} 
J. London Math. Soc. (2014) {\bf 90} (2): 436-450. 

\bibitem[K]{_Kobayashi:1976_}
Kobayashi, S., 
{\em Intrinsic distances, measures and geometric function theory,} 
Bull. Amer. Math. Soc. {\bf 82}, no. 3 (1976) 357-416. 



\bibitem[V]{_Verbitsky:ergodic_}
Verbitsky, M., 
{\em Ergodic complex structures on hyperk\"ahler manifolds,}
Acta Math. {\bf 215} (2015) 161-182.


\bibitem[Z1]{zinger1}
Zinger, A., 
{\em  Reduced genus-one Gromov-Witten invariants,} 
J. Diff. Geom. {\bf 83} (2009) 407-460. 

\bibitem[Z2]{zinger}
Zinger, A., 
{\em The reduced genus $1$ Gromov-Witten invariants of Calabi-Yau 
hypersurfaces,} 
J. A. M. S. {\bf 22} (2009) 691-737. 


\end{thebibliography}
\end{document}